\documentclass[a4paper, 11pt]{amsart}
\usepackage[T1]{fontenc}
\usepackage[latin1]{inputenc}
\usepackage[english]{babel}
\usepackage{amsmath}
\usepackage{amssymb}
\usepackage{amsfonts}

\begin{document}
\newcommand{\normi}[1]{\|#1\|}
\newcommand{\itse}[1]{\left|\,#1\right|}
\newcommand{\its}[1]{\bigl|\,#1\bigr|}
\newcommand{\rn}{\mathbb{R}^n}
\newcommand{\na}{\mathbb{N}}
\newcommand{\re}{\mathbb{R}}
\newcommand{\R}{\mathcal{R}}
\newcommand{\Z}{\mathbb{Z}}
\newcommand{\M}{\mathcal{M}}
\newcommand{\vM}{\overset{\rightarrow}{M}}
\newcommand{\ve}[1]{\overset{\rightarrow}{#1}}
\def\Xint#1{\mathchoice
 {\XXint\displaystyle\textstyle{#1}}%
 {\XXint\textstyle\scriptstyle{#1}}%
 {\XXint\scriptstyle\scriptscriptstyle{#1}}%
 {\XXint\scriptscriptstyle\scriptscriptstyle{#1}}%
 \!\int}
 \def\XXint#1#2#3{{\setbox0=\hbox{$#1{#2#3}{\int}$}
 \vcenter{\hbox{$#2#3$}}\kern-.5\wd0}}
 \def\ddashint{\Xint=}
 \def\dashint{\Xint-}
\newcommand{\m}{m}
\newcommand{\subsub}{\subset\subset}
\newcommand{\av}[3]{\underset{B(#2,#3)}{\dashint}#1(y)\,dy\,}
\newcommand{\avr}[1]{\underset{#1}{\dashint}}
\newtheorem{theorem}{Theorem}[section]
\newtheorem{lemma}[theorem]{Lemma}
\newtheorem{lause}[theorem]{Theorem}
\newtheorem{definition}[theorem]{Definition}
\newtheorem{proposition}[theorem]{Proposition}
\newtheorem{corollary}[theorem]{Corollary}
\newtheorem{question}[theorem]{Question}
\newtheorem{conjecture}[theorem]{Conjecture}
\title[Infimal Convolution in the Framework of Sobolev-functions]{On the Hamilton-Jacobi Equation and Infimal Convolution in the Framework of Sobolev-functions}
\author{Hannes Luiro}
\address{Department of Mathematics and Statistics\\
University of Jyväskylä\\
P.O.Box 35 (MaD)\\
40014 University of Jyväskylä, Finland}
\email{hannes.luiro@jyu.fi}
%\subjclass[2000]{primary 35F21, 46E35, 35B65}
%\keywords{Hamilton-Jacobi equation, Hopf-Lax formula, Sobolev spaces, regularity}
%\thanks{The author has been supported }
\maketitle
{\small \textbf{Abstract.}
We study the regularity properties of the Hamilton-Jacobi flow equation and infimal convolution in the case where initial datum function is continuous and lies in given Sobolev-space $W^{1,p}(\rn)$. We prove that under suitable assumptions it holds for solutions $w(x,t)$ that $D_xw(\cdot,t)\to Du(\cdot)$ in $L^p(\rn)$. Moreover, we construct examples showing that our results are essentially optimal.
\section{Introduction}
Consider the Hamilton-Jacobi flow equation 
\begin{equation}\label{floweq}
\begin{cases}
w_t+H(Du)=0\,\text{ in }\rn\times]0,\infty[\,,\\
w(x,0)=\lim_{t\to 0} w(x,t)=u(x)\,\text{ for every }x\in\rn\,. 
\end{cases}
\end{equation}
for initial datum function $u:\rn\to\re$ and Hamiltonian function $H:\rn\to [0,\infty[\,$.
It is known (see \cite{St}) that the unique viscosity solution $w$ of the equation (\ref{floweq}) is achieved via the Hopf-Lax formula
\begin{equation}\label{hopf}
w(x,t)=\inf_{y\in\rn}\bigg(\,u(y)+tL\big(\frac{y-x}{t}\big)\,\bigg)\,,
\end{equation}
where
\begin{equation}\label{legendre} 
L(q)=\sup_{p\in\rn}\big(\,p\cdot q-H(p)\,\big)
\end{equation}
is the Legendre-Fenchel transform of $H$. Here one of course needs to require some regularity properties for functions $u$ and $H$. It is shown in \cite{St} that (\ref{hopf}) is a solution of (\ref{floweq}) if $H$ is a strictly convex $C^1$-function with $H(0)=DH(0)=0$, and superlinear, i.e. $\frac{H(p)}{|p|}\to\infty$ as $p\to\infty$. From here on we will assume that these are valid for $H$.

For initial value function $u$ our standing assumptions are that $u$ is continuous and bounded function. In this case it is clear that $w(\cdot,t)$ provides a good approximation of $u$ when $t$ is small. Indeed, formula (\ref{hopf}) above, for given $t>0$, is also often called as an \textit{infimal convolution} of $u$ (respect to function $L$) and appears as a standard smoothing operator in the theory of viscosity solutions. In this setting the notation $w(x,t)=u_t(x)$ is typically used and so do we as well. 

We will study the properties of the infimal convolution also in the case, where $L$ needs not to be a Legendre transform of any $H$ with aforementioned properties and so that the formula (\ref{hopf}) can not be interpreted as a solution of (\ref{floweq}). Our standing assumptions for function $L$ will be that 
\begin{equation}\label{final}
L \text{ is } \text{ non-negative }C^1\text{-function},\,\, L(0)=0=DL(0) \text{ and } L \text{ is superlinear}.
\end{equation}
These assumptions almost coincide with those we made for function $H$, only the convexity assumption is dropped. It is also well known that if $L$ is the Legendre transform of $H$, then $L$ is convex function and also satisfies the required properties in (\ref{final}). 

The regularity issues for the solutions of Hamilton-Jacobi equations (or infimal convolutions) have been widely studied, see, e.g., \cite{BDR}, \cite{H}, \cite{K}, \cite{VHT} and \cite{YTW}. 
In this work we focus on the convergence problems for the derivatives of the solutions. As we know that convergence $w(\cdot, t)=u_t\to u$ holds, it is natural to ask if this convergence applies also to the derivatives. In classical setting, where $u\in C^k$, $k\in\na$ and the pointwise convergence is studied, this question is rather easy. It follows e.g. from the method of characteristics (with additional assumption $H\in C^2$) that answer is positive even for higher derivatives. 

The present paper investigates the above problem in the framework of Sobolev-functions and S. More precisely, we assume that $u$ is continuous and bounded and $u\in W^{1,p}(\rn)$ for some $1<p<\infty$ and then pose the following question: Does it hold that
\[
\tag{Q}\,\,\,\,\,Du_t\to Du\,\,\text{ in }L^p(\rn)\,\,?
\]
According to our knowledge, this problem has not been studied before.

Our first main result is the following:
\begin{theorem}\label{lptulokset}
Let $u\in W^{1,p}(\rn)\,$ be bounded, $p>n$, and suppose that there exists constant $C>0$ such that 
\begin{equation}\label{aina}
|DL(x)|\leq C\frac{L(x)}{|x|}\,.
 \end{equation}
Then  
\[\normi{Du_t-Du}_p\to 0\,\,\text{ as }t\to 0\,.\]
Moreover, $\normi{Du_t}_p\leq C(p,n,L)\normi{Du}_p$ for every $t>0$ and $u\in W^{1,p}(\rn)\,$.
\end{theorem}
It is easy to check that in the case $L(x)=C|x|^{q}$, $q>1$, appearing as the most commonly used in the theory of viscosity solutions, assumption (\ref{aina}) is valid. More careful analysis show that e.g. convexity with condition 
\begin{equation}\label{paha1}
\sup_{\frac{r}{2}<|x|<r}L\leq C\inf_{\frac{r}{2}<|x|<r}L\,
\end{equation}
for all $r>0$ and constant $C<\infty$ is sufficient for (\ref{aina}).

Regarding the solutions of the Hamilton-Jacobi equation, Theorem \ref{lptulokset}, combined with Proposition \ref{pakko} in section \ref{osa2}, gives the following Corollary:
\begin{corollary}\label{seuraus1}
Let $p>n$, $u\in W^{1,p}(\rn)$ bounded and assume that $H$ satisfies our standing assumptions (given above). Additionally, suppose that there exists $C>1$ and $C'>1$ such that 
\begin{equation}\label{paha2a}
\frac{H(2x)}{H(x)}>\,2C\text{ for every }x\in\rn\,\text{ and }
\end{equation} 
\begin{equation}\label{paha3a}
\sup_{\partial B(0,r)}\,H\,< C'\inf_{\partial B(0,r)}\,H\,.
\end{equation}
for every $r>0$. Then $D_xw(\cdot,t)\to Du(\cdot)$ in $W^{1,p}(\rn)$ if $t\to 0\,$, and $\normi{D_xw(\cdot,t)}_p\leq C(p,n,H)\normi{Du}_p$ for every $t>0$ and $u\in W^{1,p}(\rn)\,$.
\end{corollary}
\noindent The proofs of the above results are given in section \ref{osa2}.

In section \ref{osa3} we show that our assumptions in Theorem \ref{lptulokset} (and Corollary \ref{seuraus1}) are essentially sharp. 
%It is clear that in the case of inf-convolution, i.e. $L(x)=C|x|^2\,$ or more generally, when $L(x)=C|x|^{q}$, $q>1$ %assumption (\ref{aina}) is valid. More careful analysis show that e.g. convexity with condition 
%\begin{equation}\label{paha1}
%\inf_{\frac{r}{2}<|x|<r)}L\leq C\sup_{\frac{r}{2}<|x|<r}L\,
%\end{equation}
%for all $r>0$ and constant $C<\infty$ is sufficient for (\ref{aina}). Furthermore, we will verify that
%\begin{proposition}
%If $L$ is a Legendre transform of $H$ (with aforementioned properties) then (\ref{aina}) follows if there exists $C>1$ so %that
%\begin{equation}\label{paha2}
%\frac{H(2x)}{H(x)}>\,2C\text{ for every }x\in\rn\,
%\end{equation}
%and moreover, there exists $C'>1$ so that 
%\begin{equation}\label{paha3}
%\inf{B(x,C'r)}^c\, H>\,\sup_{B(x,r)}\,H\,,
%\end{equation}
%for every $r>0$.
%\end{proposition}
Relating to the assumption $p>n$, it is shown (even in the most important special case $L(q)=\frac{|q|^2}{2}$) that 
\begin{theorem}\label{ekavasta}
If $n\geq 2$, then there exists continuous and compactly supported function $u\in W^{1,n}(\rn)$ such that $|Du_t|\to \infty$ when $t\to 0$ in a set of positive measure.
\end{theorem} 
\noindent This implies that answer to (Q) is definitely negative when $1\leq p\leq n\,$. 

The rest of the section \ref{osa3} deals with assumption (\ref{aina}).
The most natural case where (\ref{aina}) fails appears when $L$ behaves in infinity like exponential function, which corresponds to the case where $H$ is of type $|x|\log(|x|)\,$ and (\ref{paha2a}) does not hold for $H$. In the second example of section three it is shown that in this case the convergence (Q) fails in general also in the case $p>n\geq 3$. In this connection we also state some open questions, especially relating to the cases $n=1$ and $n=2$.

We end up with justifying in our third example the assumption of 'quasi-radiality', i.e. property (\ref{paha3a}) by considering the example in $\re ²$, where $H(x,y)=C_s|x|^{s}+C_{s'}|y|^{s'}$, $s>s'>1$. This appears as a model case for the failure of (\ref{paha3a})
%, as well as for the failure of (\ref{paha1}) for the Legendre transform $L$ of $H$ 
($H$ still satisfies other required properties). We show that convergence (Q) does not hold in this situation for any $1<p<\infty$.\\

\textit{Acknowledgements.}
The author would like to thank Peter Lindqvist for pointing me out the open problems studied in this paper as well as Vesa Julin, Petri Juutinen and Eero Saksman for useful discussions.

\section{Boundedness and convergence of the gradients}\label{osa2}

Before proving Theorem \ref{lptulokset}, we have to verify some auxiliary lemmas. Althought most of these results are rather easy and well known in optimization theory, we have not find sufficiently convenient references for our purposes. Because of that we prefer to express the complete proofs. 

Firstly, it is well known that for every fixed $t>0$ functions $u_t$ are locally Lipschitz and in our case where $u$ is assumed to be bounded, functions $u_t$ are even globally Lipschitz. More precisely, suppose that $C=\sup u$ and choose $R>0$ such that $tL(a)\geq 2 C\,$ if $|a|\geq R\,$. Then one can easily verify that $u_t$ is Lipschitz with constant
\[\sup_{a\in B(0,R)}|DL(a)|\,.\] 
Here continuity of function $u$ does not have any role.

Next we remark that the Hopf-Lax formula may be written in many different ways. In this work we will exploit the following equivalent definition:
\begin{equation}
u_t(x)=\inf_{a\in\rn}\,\big( u(x+ta)+tL(a)\,\big)\,.
\end{equation}

Furthermore, since we assumed that function $u$ is bounded and $L(x)\to\infty$ as $|x|\to\infty$, it follows that for every $t>0$ there exists $R_t>0$ such that
\begin{equation}\label{lokaalisuus}
u_t(x)=\inf_{a\in B(0,R_t)}\,\big(u(x+ta)+tL(a)\,\big)\,.
\end{equation}
for every $x\in \rn$. Assuming additional properties for $L$ or $u$ would make it possible to verify more precise estimates for $R_t$. 

 %Remark that continuity of $u$ does not play any role here. 
%From this one can quite easily prove that functions $u_t$ are globally Lipschitz with constant $C_k\leq \sup_{a\in %B(0,R_t)}|DL()|$ where Let us begin with verifying that functions $u_t$ are locally Lipschitz-functions for every $t>0$. %More precisely, suppose that $x,z\in\rn$ and $0<|x-z|<1$. Suppose that $u_t(x)>u_t(z)$ and sequence $y_k$ gives the value $u_t(z)$. Then we simply estimate the difference $u_t(x)-u_t(z)$ by choosing the same sequence at point $z$, thus
%\begin{align*}
%u_t(x)-u_t(z)&=u_t(x)-\big(\lim_{k\to\infty}\, u(y_k)+tL\big(\frac{z-y_k}{t}\big)\,\big)\\&\leq 
%\lim_{k\to\infty}\,\big(\, u(y_k)+tL\big(\frac{x-y_k}{t}\big)\,\big)
%-\lim_{k\to\infty}\,\big(\, u(y_k)+tL\big(\frac{z-y_k}{t}\big)\,\big)\\
%&=t\,\lim_{k\to\infty}\big(L\big(\frac{x-y_k}{t}\big)-L\big(\frac{z-y_k}{t}\big)\big)\,.
%\end{align*}
%Then it follows from (\ref{lokaalisuus}) that $\frac{|z-y_k|}{t}\leq\frac{R}{t}$ when $k$ is big enough, thus %$\frac{|x-y_k|}{t}\leq \frac{R}{t}+\frac{1}{t}$ and we get from above that   
%\begin{align*}
% u_t(x)-u_t(z)&\leq \lim_{k\to\infty}\, %t\,\big|\frac{x-y_k}{t}-\frac{z-y_k}{t}\big|\bigg(\,\sup_{B(0,\frac{R}{t}+\frac{1}{t})} |DL|\,\bigg)\\
%&=|x-z|\bigg(\,\sup_{B(0,\frac{R}{t}+\frac{1}{t})} |DL|\,\bigg)\,.
%\end{align*} 

Throughout this work the following notation is used:
For every $x\in\rn$ define set $\R u_{t}(x)$ such that 
$a\in \R u_{t}(x)$ if 
%and only if there exists sequence of points $a_k$ converging to $a$ so that
%\[
%u_{t}(x)=\lim_{k\to\infty}\,u(x+a_kt)+tL(a_k)\big)\,.
%\]
%It is easy to check that $R u_t(x)\neq \emptyset$ (that follows for example from (\ref{lokaalisuus})). Moreover
%If $u$ is continuous then it also clearly holds that 
\begin{equation}\label{loppu2}
u_{t}(x)=u(x+at)+tL(a)\,.%\text{ for every }a\in\R u_{t}(x)\,
\end{equation}
It is easy to check that $\R u_t(x)$ is always non-empty and compact set.
%but for general $u$ this of course may not hold.
The following inequality (suppose $|u|\leq C$) 
\begin{equation}\label{loppu1}
 u_t(x)=u(x+ta_t)+tL(a_t)\geq -C+t|a_t|\frac{L(a_t)}{|a_t|},
\end{equation}
if $a_t\in \R u_t(x)$, combined with the superlinearity of $L$ (and the fact that $a_t$ minimizes (\ref{loppu2})), guarantees that for any given $x$ the sets $t\R u_t(x)$, $t>0$, are uniformly bounded (respect to $t$). This, combined with the continuity of $u$, also implies the (pointwise) convergence $u_t\to u$ as $t\to 0$.

Since functions $u_t$ are Lipschitz, they are differentiable almost everywhere by Rademacher's theorem. Actually more careful analysis show that they are even twice differentiable a.e. This is due to the fact that one can show that functions $u_t$ are semiconcave (see e.g. \cite{JJ}).

It is easy to see that functions $u_t$ do not need to be differentiable everywhere. However, the directional derivatives exist at every point, which is verified in the following lemma.  

\begin{lemma}\label{directional}
Suppose that $u$ is continuous and bounded function and $t>0$. Then $u_t$ has directional derivatives at every point to any direction $\gamma$ and
\begin{equation}
D_{\gamma}u_t(x)=\inf_{a'\in\R u_t(x)}(-D_{\gamma}L(a'))=-D_{\gamma}L(a)\,\textit{ for some }a\in\R u_t(x)\,.
\end{equation}
\end{lemma}
\textit{Proof.} Let $t>0$ and $x\in\rn$. Recall first that by compactness of set $\R u_t(x)$ there exists $a\in\R u_t(x)$ such that 
\begin{equation}
-D_{\gamma}L(a)=\inf_{a'\in\R u_t(x)}(-D_{\gamma}L(a'))\,.
\end{equation}
Moreover, continuity of $u$ implies that $u_t(x)=u(x+ta)+tL(a)$.
We first estimate the difference $u_t(x+h\gamma)-u_t(x)$, $h>0$,  from above by choosing points $a_h$ such that $a_h=a-\frac{h\gamma}{t}$. It follows that
\begin{align*}
\frac{u_t(x+h\gamma)-u_t(x)}{h}&\leq
\frac{1}{h}\big(\,(u(x+h\gamma +ta_h)+tL(a_h))-(u(x+ta)+tL(a))\,\big)\\
&=\frac{t}{h}\big(\,L(a_h)-L(a)\,\big)=\frac{\,L(a-\frac{h\gamma}{t})-L(a)}{\frac{h}{t}}\,.
\end{align*}
Taking the limit yields that
\begin{equation}
\limsup_{h\to 0^+}\,\frac{u_t(x+h\gamma)-u_t(x)}{h}\leq -D_{\gamma}L(a)\,.
\end{equation}

Assume then that there exists sequence $h_k>0$, $h_k\to 0$ as $k\to\infty$, and $\lambda >0$ so that  
\begin{equation}\label{AT}
\frac{u_t(x+h_k\gamma)-u_t(x)}{h_k}< -D_{\gamma}L(a)-\lambda\, 
\end{equation}
for every $k$. Suppose that $a_k\in \R u_t(x+h_k)$ and choose $a'_k=a_k+\frac{h_k\gamma}{t}\,$. Recall that sequence ($a_k$) has to be bounded and by extracting a subsequence, if needed, we may assume that $a_k\to a_0$. Moreover, 
\[
u_t(x)=\lim_{k\to\infty} u_t(x+h_k)=\lim_{k\to\infty}\,u(x+h_k+ta_k)+tL(a_k)=\lim_{k\to\infty}\,u(x+ta'_k)+tL(a'_k)
\]
which guarantees that $a_0\in\R u_t(x)\,$. Then we estimate the difference from below using points $a'_k$: 
\begin{align*}
\frac{u_t(x+h_k\gamma)-u_t(x)}{h_k}&\geq
\frac{1}{h_k}\big(\,(u(x+h_k\gamma +ta_k)+tL(a_k))-(u(x+ta'_k)+tL(a'_k))\,\big)\\
&=\frac{t}{h_k}\big(\,L(a_k)-L(a'_k)\,\big)=\frac{\big(\,L(a_k)-L(a_k+\frac{h_k\gamma}{t})\,\big)}{\frac{h_k}{t}}\,. 
\end{align*}
Then recall that $L$ is $C^1$-function and $a_k\to a_0$ as $k\to\infty$. This implies that above the last quantity converges to $-D_{\gamma}L(a_0)$. Because $a_0\in \R u_t(x)$ we get that
\begin{align*}
\liminf_{k\to\infty}\frac{u_t(x+h_k\gamma)-u_t(x)}{h_k}\geq -D_{\gamma}L(a_0)\geq -D_{\gamma}L(a).
\end{align*}
This contradicts with (\ref{AT}) and the proof is complete.
\hfill$\Box$\\ 
 
We get the following corollary in the points of differentiability of $u_t\,$:
\begin{corollary}\label{peruskaava}
Let $u$ be continuous and bounded and $t>0$. Then $u_t$ is differentiable at $x$ and $Du_t(x)=-DL(a)$ for every $a\in \R u_t(x)$ if and only if  
$DL(a)=DL(a')$ for every $a,a'\in\R u_t(x)\,$.
\end{corollary}
\textit{Proof.}
Suppose, on the countrary that $u_t$ is differentiable at $x$ such that there exists $a_1,a_2\in \R u_t(x)$ so that $DL(a_1)\neq DL(a_2)\,$. Then there exists direction $\gamma$ so that $DL(a_1)\cdot\gamma < DL(a_2)\cdot\gamma\,$. Then we get by Lemma (\ref{directional})  that $D_{\gamma}u_t(x)\leq -DL(a_2)\cdot\gamma<-DL(a_1)\cdot\gamma=DL(a_1)\cdot(-\gamma)\leq-D_{-\gamma}u_t(x)\,=D_{\gamma}u_t(x)\,$, which is a contradiction. Another direction follows by the same argument. 
\hfill$\Box\,$\\

\textit{Remark 1.} Though we know that functions $u_t$ are differentiable almost everywhere, we prefer writing our lemmas or theorems using directional derivatives $D_{\gamma}u_t$ since they really exist at every point for \textit{every} $t\in(0,1)$. Then we can avoid possible technical problems which may occur if there were some exceptional sets of measure zero for every $t\in(0,1)\,$. 

\textit{Remark 2.}
It is often preferable to study the regularity of the solutions $u_t$ in terms of subdifferentials (see e.g. \cite{PR}). 
In this presentation we have chosen to avoid this since they seem not to give any advantage for our short-term goals. 

\textit{Remark 3.} In above results it was assumed that $u$ is continuous function. However, this is actually not needed but simplifies the proof since it makes possible to write $u_t(x)=u(x+ta)+tL(a)$ if $a\in\R u_t(x)$ instead of using the sequences $(a_k)$ giving the infimum. 

\textit{Remark 4.} It is clear that also the boundedness of $u$ can be replaced by assuming only that the growth of $u$ in infinity is suitably comparable to the growth of $L$. This remark (as well as the remarks above) is involved in forthcoming results, as well.\\

To prove the convergence $Du_t\to Du$ in $L^p$ when $p>n$ we first show that this convergence holds pointwise a.e. This part follows basically from the next theorem.

\begin{theorem}\label{classical}
Suppose that $u$ is continuous and bounded function and differentiable at $x\in\rn\,$. Then $D_{\gamma}u_t(x)\to Du(x)\cdot\gamma$ as $t\to 0$ for any direction $\gamma$. 
\end{theorem}
\textit{Proof.} Denote that $Du(x)=:D$. By the differentiability of $u$ at $x$,  
there exists function $\varepsilon_x$ such that $\varepsilon_x(h)\to 0$ as $|h|\to 0$ and 
\begin{equation}\label{muoto2}
u(x+ta)=u(x)+t(D\cdot a)+t\varepsilon_x(ta)|a|\,.
\end{equation}
Suppose that $a_t\in\R u_t(x)$ so that $D_{\gamma}u_t(x)= -D_{\gamma}L(a_t)$. To complete the proof we have to show that $D_{\gamma}L(a_t)\to -D\cdot \gamma$. 
It follows from (\ref{muoto2}) that for every $t>0$, point $a_t$ minimizes the quantity
\begin{equation}\label{muoto}
D\cdot a+\varepsilon_x(ta)|a|+L(a)\,.
\end{equation}
As it was deduced from (\ref{loppu1}) above, $\{ta_t\}_{t>0}$ is bounded. Combining this with the fact that $a_t$ minimizes (\ref{muoto}) and using again the superlinearity of $L$ one obtains that even $\{a_t\}_{t>0}$ is bounded. 

Next we denote by $A$ the set of points $a'$ which are the minimizers of the above quantity without the error term $\varepsilon_x(ta)|a|$, thus every $a'\in A$ minimizes
\begin{equation}
D\cdot a+L(a).
\end{equation} 
Superlinearity of $L$ implies that $A\not=\emptyset$.
Then it clearly holds (recall that $L$ is $C^1$-function) that $DL(a')=-D$ for every $a'\in A$. 

Suppose then that our claim is not true, thus there exists a sequence $a_{t_k}$ such that $D_{\gamma}L(a_{t_k})\not\to -D\cdot \gamma$. Since the sequence $a_{t_k}$ has to be bounded we may assume that it converges, say to $a_0'$. Then 
\begin{align*}
D\cdot a_0'+L(a_0')&=\lim_{k\to\infty}\big(D\cdot a_{t_k}+\varepsilon_x(t_ka_{t_k})|a_{t_k}|+L(a_{t_k})\big)\\
&\leq \lim_{k\to\infty}\big(D\cdot a+\varepsilon_x(t_ka)|a|+L(a)\big)=D\cdot a+L(a)\,
\end{align*}
for every $a\in\rn$, thus $a'_0\in A$ and $DL(a'_0)=-D$. Therefore 
\[D_{\gamma}L(a_{t_k})\to D_{\gamma}L(a_0')=-D\cdot\gamma\,.\] 
This is the desired contradiction.
\hfill$\Box$.\\

\textit{Remark 1.} 
This implies that if $u$ is differentiable almost everywhere, like if $u\in W^{1,p}(\rn)$, $p>n\,$, then $D_{\gamma}u_t\to Du\cdot\gamma$ almost everywhere and $Du_{t_k}\to Du$ a.e for any sequence $t_k\to 0\,$ (here it is reasonable to formulate the latter claim for sequences, see the remarks above).

\textit{Remark 2.} Above theorem accounts for our assumption $L\in C^1\,$ (or strictly convexity of $H$). Dropping this assumption would yield problems with this theorem. The main reason for this is the fact that we really need the above theorem in the form where anything else about the smoothness of $u$ around point $x$ is not known than the differentiability in that single point $x\,$.\\

The remaining ingredient of our main result relies to the following rather standard inequality (following e.g. from Poincare inequality, see \cite{A}).
\begin{lemma}\label{perus2}
Suppose that continuous function $u:\rn\to\re$ has locally integrable weak partial derivatives and $p>n\,$. Then 
\[|u(x)-u(y)|\leq\,C(p,n)|x-y|\,\bigg(\avr{B(x,2|x-y|)}|Du|^p\bigg)^{\frac{1}{p}}\]
for every $x,y\in\rn$ such that $B(x,2|x-y|)\subset\Omega$.
\end{lemma}

Using the above inequality we get the desired $p$-integrable majorant for functions $Du_t\,$ when $|DL(x)|\lesssim \frac{L(x)}{|x|}\,$. Below $M$ denotes the classical centered Hardy-Littlewood maximal function.
\begin{theorem}\label{bound}
If $u$ is bounded and continuous and $L$ satisfies $|DL(x)|\leq C\frac{L(x)}{|x|}\,$, then 
\begin{equation}
\sup_{\gamma}|D_{\gamma}u_t(x)|\leq C'\big(M(|Du|^p)(x)\big)^{\frac{1}{p}}\,\text{ for every }t\in]0,1]\,\text{ and }x\in\rn\,.
\end{equation}
\end{theorem}
\textit{Proof.} 
Observe first that by Lemma (\ref{directional}) 
\[\sup_{\gamma}|D_{\gamma}u_t(x)|=\sup_{\gamma}|D_{\gamma}L(a_{\gamma})|=|D_{\gamma}L(a)|\,,\]
for some $a\in\R u_t(x)\,$ (here we also used the fact that $\R u_t(x)$ is compact and $L$ is $C^1$-function). 
Then, if $a=0$, the claim is trivially true since $DL(0)=0\,$.
Suppose then that $a\neq 0$  whence it holds that 
\begin{align*}
u(x+ta)+ tL(a)\leq u(x)\,\,
\Longrightarrow\,\,& \frac{L(a)}{|a|}\leq\frac{|u(x+ta)-u(x)|}{t|a|}\,.
\end{align*}
Then we get by Lemma \ref{perus2} and assumption $|DL(x)|\leq C\frac{L(x)}{|x|}\,$ that
\begin{align*}
|DL(a)|\leq C\frac{L(a)}{|a|} \leq 
&\,CC(p,n)\bigg(\avr{B(x,2|x-y|)}|Du|^p\bigg)^{\frac{1}{p}}\leq C'\big(M(|Du|^p)(x)\big)^{\frac{1}{p}}\,
\end{align*}
for every $x$.
\hfill$\Box$\\

Then we are ready to prove Theorem \ref{lptulokset}.

%\begin{theorem}
%If $p>n$ and $L$ satisfies $|DL(x)|\leq C\frac{L(x)}{|x|}\,$, then 
%\begin{equation}
%\normi{Du_t}_p\leq C'\normi{Du}_p\,.
%\end{equation}
%\end{theorem}
%\textit{Proof.}
%Let us choose $c=\frac{1+\frac{p}{n}}{2}>1\,$ whence it holds that $\frac{p}{c}>n$. Then combining the boundedness of the %maximal operator in $L^c(\rn)$ and Theorem \ref{perus3} we get that
%\begin{align*}
%\normi{Du_t}_p &\leq 2C(n,p)\normi{\big(M(|Du|^{\frac{p}{c}})\big)^{\frac{c}{p}}}_p
%\,=\,2C(n,p)\big(\normi{M(|Du|^{\frac{p}{c}})}_c\big)^{\frac{c}{p}}\\
%&\leq C'(n,p)\big(\normi{|Du|^{\frac{p}{c}}}_c\big)^{\frac{c}{p}}\,=\,C'(n,p)\normi{|Du|}_p\,.
%\end{align*}

\subsection*{Proof of Theorem \ref{lptulokset}}
Let $p>n$, $u\in W^{1,p}(\rn)$ bounded, and suppose that $L$ satisfies $|DL(x)|\leq C\frac{L(x)}{|x|}\,$. We have to show that  
\[\normi{Du_t-Du}_p\to 0\,\,\text{ as }t\to 0\,,\]
and, moreover, that $\normi{Du_t}_p\leq C'\normi{Du}_p$ for every $t>0$.

Let us begin with choosing a sequence $t_k>0$ and $t_k\to 0$ as $k\to\infty$.
Recall that in case $p>n$ function $u$ is a.e. differentiable. Then it follows by Theorem (\ref{classical}) that $Du_{t_k}\to Du$ almost everywhere. Then the first of the claims above follows by Lebesgue dominated convergence theorem if we only find the desired $p$-integrable bound for the sequence $|Du_{t_k}|$. For that, let us choose $c=\frac{1+\frac{p}{n}}{2}>1$ and apply Theorem (\ref{bound}) for exponent $\frac{p}{c}$ ( where $n<\frac{p}{c}<p$ ) to obtain that   
\[|Du_{t_k}|\leq CM(|Du|^{\frac{p}{c}})\big)^{\frac{c}{p}}\,\,\text{ a.e}.\]
Furthermore, using the boundedness of the maximal operator in $L^c(\rn)$ we get  
\begin{align*}
\normi{\big(M(|Du|^{\frac{p}{c}})\big)^{\frac{c}{p}}}_p
\,=\big(\normi{M(|Du|^{\frac{p}{c}})}_c\big)^{\frac{c}{p}}
\leq C'\big(\normi{|Du|^{\frac{p}{c}}}_c\big)^{\frac{c}{p}}\,=\,C'\normi{|Du|}_p\,.
\end{align*}
These formulas also immediately imply the latter claim.
\hfill$\Box\,$\\

We end up this section with the proof of the proposition \ref{pakko}. When combined with Theorem \ref{lptulokset}, it immediately implies Corollary \ref{seuraus1}. 
\begin{proposition}\label{pakko}
If $L$ is a Legendre transform of $H$ (with aforementioned properties) then (\ref{aina}) follows if there exists $C>1$ and $C'>1$ so that
\begin{equation}\label{paha2}
\frac{H(2x)}{H(x)}>\,2C\text{ for every }x\in\rn\,\text{ and }
\end{equation} 
\begin{equation}\label{paha3}
\sup_{\partial B(0,r)}\,H\,< C'\inf_{\partial B(0,r)}\,H\,.
\end{equation}
for every $r>0$.
\end{proposition}

\textit{Proof.}
Suppose that $H$ satisfies the requirements of the proposition. It is easy to check that for every $q\in\rn$ there is unique $p_q$ such that 
\begin{equation}
L(q)=\sup_{p\in\rn} p\cdot q -H(p)\,=\,p_q\cdot q-H(p_q)\,\,\text{ and, moreover, } DL(q)=p_q\,.
\end{equation}
Thus, the claim of the proposition is equivalent with 
\begin{equation}\label{hassu1}
|p_q||q|\leq C''(p_q\cdot q-H(p_q))\,,
\end{equation}
for some constant $C''<\infty$.
Let then $p'_q:=\frac{(p_q\cdot q)q}{|q|^2}$ denote the projection of $p_q$ on a line $\{tq:t\in\re\}$. Because $p_q\cdot q=p_q'\cdot q$ , it follows that $H(p_q)\leq H(p_q')$. Using the convexity of $H$ and assumption (\ref{paha3}) we then conclude that $|p_q|\leq C'|p_q'|\,$. This in turn implies that 
\begin{equation}\label{hassu2}
|p_q||q|\leq C'(p_q\cdot q)\,.
\end{equation}
Furthermore, using (\ref{paha2}) we get that 
\begin{align*}
p_q\cdot q-H(p_q)\geq& \,\frac{p_q}{2}\cdot q-H(\frac{p_q}{2})\geq \frac{p_q}{2}\cdot q-\frac{1}{2C}H(p_q)\\
=&\,\frac{1}{2C}(p_q\cdot q -H(p_q))+\frac{C-1}{2C}(p_q\cdot q)\,.
\end{align*}
Finally, we conclude with 
\begin{equation}\label{hassu3}
p_q\cdot q\leq \frac{2C-1}{C-1}(p_q\cdot q-H(p_q))\,. 
\end{equation}
This completes the proof, since (\ref{hassu1}) follows from (\ref{hassu2}) and (\ref{hassu3}).
\hfill$\Box$\\

This proposition shows that assumptions (\ref{paha2}) and (\ref{paha3}) for $H$ guarantee that Legendre transformation $L$ of $H$ satisfies (\ref{aina}) and Theorem \ref{lptulokset} is usable. 
    
\section{Counterexamples}\label{osa3}
In this section we consider the sharpness of the results we proved in the previous section. This section is splitted in three subsections. The first one deals with assumption $p>n$ while two latter subsections are devoted to the assumption $DL(x)\lesssim \frac{L(x)}{|x|}\,$. 
\subsection*{The case $p\leq n$}
The following example deals with case $p\leq n$ and the most important special case $L(q)=\frac{1}{2}|q|^2=H(q)$. If we only want to show that theorem (\ref{lptulokset}) does not hold in this case, thus the convergence or the boundedness of the gradients $Du_t$ in Sobolev-norm does not hold, we may do this much more simple than the following construction by using suitable radial functions. However, in the following example it is shown that in this case even pointwise convergence may totally fail. 

For this example we first recall that in the case $L(q)=\frac{1}{2}|q^2|$ it follows directly from Corollary \ref{peruskaava} that for a.e. $x\in\rn$ 
\begin{equation}\label{kaavakaava}
Du_t(x)=\frac{x-y^*}{t}\,\,
\end{equation}
for every $y^*$ such that the infimum in (\ref{hopf}) is attained at $y^*$. 

In the core of the construction below is (of course) the fact that functions $u\in W^{1,n}(\rn)$, $n\geq 2$, may be even discontinuous , i.e. they do not need to have continuous representatives. Therefore it is easy to construct a continuous function $u\in W^{1,n}(\rn)$ for which 
\[\liminf_{k\to\infty}\bigg(\sup\big{\{}\frac{u(x)-u(y)}{|x-y|}\,:\,y\in B(x,2^{-k})\setminus B(x,2^{-k-1})\big{\}}\bigg)=\,\infty\,\] 
even for almost every $x\,$. On the other hand, the careful analysis of the definition of $u_t$ suggests that $\liminf_{t\to 0}|Du_t(x)|$ should be typically comparable to this quantity. In the following example this intuition
is verified.

\subsubsection*{Proof of Theorem \ref{ekavasta}}
Let $n\geq 2$ and $L(q)=\frac{1}{2}|q|^2$. We have to find continuous and compactly supported function $u\in W^{1,n}(\rn)$ such that $|Du_t|\to \infty$ in a set of positive measure when $t\to 0$.
Let us begin with defining grid sets $A_k$ in $[0,1]^n$ by 
\begin{equation}
A_k=\bigg(\{\frac{i}{4^k}: 0\leq i\leq 4^k\,,\,i\in\mathbb{Z}\}\bigg)^n\,.
\end{equation}
Then denote by $N_k$ the number of points in $A_k$. Define also radially increasing functions $g^k$ such that $g^k(0)=-\frac{1}{2^k}$,
$\normi{Dg^k}_n\leq\frac{1}{2^kN_k}\,$ and $g^k$ is supported in $B(0,r_k)$. Moreover, radii $r_k$ are chosen to be so small that $r_k<<\frac{1}{4^k}$ and, especially,  
\[\bigg|\bigcup_{k=1}^{\infty}\bigcup_{a\in A_k}B(a,2r_k)\bigg|=:|B|<\frac{1}{2}.\] 
Finally, define 
\begin{equation}
u(x)=\inf_{k\in\na}\bigg(\sum_{a\in A_k}g^k(x-a)\bigg)=:\inf_{k\in\na}f^k(x)\,.
\end{equation}
Observe that the distance between two points in $A_k$ is always greater than $r_k$ which implies that for every $x$ and $k$ there exists at most one $a\in A_k$ such that $g^k(x-a)\neq 0$. Moreover, $u$ is continuous function since the functions $f^k$ are continuous and $|f^k|\leq 2^{-k}$ for every $k$. Furthermore, it is easy to check that  
\begin{align*}
\normi{Du}_{n}\leq \sum_{k=1}^{\infty}N_k\normi{Dg^k}_n\leq \sum_{k=1}^{\infty}N_k\frac{1}{2^kN_k}\,=1\,.
\end{align*}
This implies that $u\in W^{1,n}(\rn)$.

In what follows we prove that $Du_t\to\infty$ a.e in $B^c$.
For the proof, observe first that
\begin{align*}
u_t(x)&=\inf_{y\in\rn}\,\big(\,u(y)+\frac{|x-y|^2}{2t}\,\big)=
\inf_{y\in\rn}\big(\inf_{k\in\na}\big(f^k(y)+\frac{|x-y|^2}{2t}\big)\,\big)\\
&=\inf_{k\in\na}\big(\inf_{y\in\rn}\big(f^k(y)+\frac{|x-y|^2}{2t}\big)\,\big)\,=\,\inf_{k\in\na}f^k_t(x)\,
\end{align*}
for every $x\in\rn$, $t>0$.

Suppose then that $t\in[4^{-j},4^{-j+1}]$ and $k>\frac{3}{4}j=:k_0$. Then there exists point $y^*\in A_{k_0}$ so that $|x-y^*|\leq C_n4^{-k_0}$ whence 
\[\frac{|x-y^*|^2}{2t}\leq C'_n\frac{4^{-2k_0}}{4^{-j}}< C'_n\frac{1}{2^j}.\]
This implies that
\begin{align*}
f^k_t(x)&\geq -\frac{1}{2^k}\geq -\frac{1}{2^{k_0}}+\frac{1}{2^{k_0+1}}=-\frac{1}{2^{k_0}}+\frac{2^{\frac{j}{4}-1}}{2^j}\\ &\geq -\frac{1}{2^{k_0}}+C'_n\frac{1}{2^j}
\geq f^{k_0}(y^*)+\frac{|x-y^*|^2}{2t}\geq f_t^{k_0}(x).
\end{align*}

We conclude that if $t\in[4^{-j},4^{-j+1}]$ then there exists integer $k_0$ such that 
$u_t(x)=f^{k_0}_t(x)$ and $k_0\leq \frac{3}{4}j\,$. Especially, if $u_t(x)$ is achieved in point $y^*$ then it follows from above that $f^k(y^*)\neq 0$ for some $k\leq\frac{3}{4}j$, thus
\begin{equation}\label{helppo}
y^*\in\bigcup_{k=1}^{\frac{3}{4}j}\bigg(\bigcup_{a\in A_k}B(a,r_k)\bigg)\,.
\end{equation}
Suppose then that $u_t(x)$ is achieved in $y^*$ such that 
$|x-y^*|\leq 4^{-\frac{5}{6}j}$ whence it holds that 
  
\begin{align*}
x&\in\bigcup_{k=1}^{\frac{3}{4}j}\bigg(\bigcup_{a\in A_k}B(a,r_k+4^{-\frac{5}{6}j})\bigg)\\
&\subset\bigg(\bigcup_{k=1}^{\infty}\bigcup_{a\in A_k}B(a,2r_k)\bigg)\,\bigcup\,
\bigg(\bigcup_{k=1}^{\frac{3}{4}j}\bigcup_{a\in A_k}B(a,2(4^{-\frac{5}{6}j}))\bigg)\\
&= B\cup \bigcup_{a\in A_{\frac{3}{4}j}}B(a,2(4^{-\frac{5}{6}j}))=:B\cup B_j\,.
\end{align*}
Observe then that if $x\in B^c\cap B_j^c$ then $|x-y^*|>4^{-\frac{5}{6}j}$ and (recall (\ref{kaavakaava}))  
\[|Du_t(x)|\overset{a.e.}{=}\frac{|x-y^*|}{t}>4^{\frac{1}{6}j}\] 
if $t\in[4^{-j}, 4^{-j+1}]\,$. It follows that
 \[\{x\in B^c: Du_t(x)\overset{t\to 0}{\not\to}\infty\}\subset\bigcap_{m=1}^{\infty}\bigcup_{j=m}^{\infty}B_j\,.\]
Finally, observe that $|B_j|\leq C_n(4^{\frac{3}{4}j}4^{-\frac{5}{6}j})^n\leq C_n 4^{-\frac{1}{12}j}$ and so  $\sum_{j=1}^{\infty}|B_j|<\infty$. This guarantees that $Du_t\to\infty$ a.e. in $B^c$.
\hfill$\Box$\\

\subsection*{Exponential but radial growth}
As it was mentioned in the introduction, perhaps the most natural case when assumption 
$|DL(x)|\lesssim\frac{|L(x)|}{|x|}$ does not hold appears when $L$ behaves in infinity like exponential function. 
At first glance it may look that the more faster $L$ grows in infinity the better the local regularity of $u$ is retained for function $u_t$. In certain sense this is true, especially at the points where $u$ is differentiable. However, in case of Sobolev-functions, even if $p>n$, there might exist singular points for the gradient like when $u(x)=|x|^{\alpha}$, where $\alpha>1-\frac{n}{p}$. It turns out that around these singular points the bad behaviour of $L$ at infinity may kill the uniform local $p$-integrability of the derivatives $Du_t$. More precisely, we prove the following theorem: 

\begin{theorem}
Let $p>n\geq 3$ and suppose that $C^1$-function $L:\rn\to[0,\infty]$ coincides with function $e^{|x|}$ when $|x|\geq 1\,$. In $B(0,1)$ suppose that $L$ is strictly convex with $L(0)=DL(0)=0\,$.     
Suppose that functions $u_t$ are determined by $L$ as in (\ref{hopf}).
Then there exists compactly supported function $u\in W^{1,p}(\rn)$ such that $\normi{Du_t}_{p}\to\infty$ as $t\to 0\,$.
\end{theorem} 
\textit{Proof.}
Let us fix $\frac{1}{2}<\alpha<1$ such that $(\alpha-1)p+n>0\,$. Moreover, define for every $k\in\na$ constants $C_k$ by 
\[C_k=\frac{2^{k(\alpha-1+\frac{n}{p})}}{k^{\frac{3}{2p}}}\,.\]
and functions $u^k$ in $\re$ by 
\begin{equation}\label{definit}
u^k(x)=\begin{cases}
        &C_k|x|^{\alpha}\,\, \text{ if } |x|<2^{-k}\,.\\
	& C_k2^{(-k)\alpha}\,\,\text{ otherwise}\,.
       \end{cases}
\end{equation}
Let us then assume that $t\in[2^{-k},2^{-k+1}]$ and consider the function $u^k_t$,
\begin{equation}\label{definit2}
 u^k_t(x)=\inf_{y\in\re}(u(y)+te^{\frac{|x-y|}{t}})\,.
\end{equation}
We are going to show that actually $u^k_t$ coincides with $te^{\frac{|x|}{t}}$ in $B(0,r)$ such that 
\[r\geq t\log(\frac{t^{\alpha-1}C_k}{2})=:r_0\,.\] 
To prove this, suppose that $|x|=r_0\,$. Then it follows that 
\begin{align*}
|(te^{\frac{\cdot}{t}})'(x)|=e^{\frac{|x|}{t}}= e^{\frac{r_0}{t}}=\frac{t^{\alpha-1}C_k}{2}\leq \frac{(2^{-k})^{\alpha-1}C_k}{2}\,\leq \alpha(2^{-k})^{\alpha-1}C_k\,.
\end{align*}
Since we know that the derivative of $tL\big(\frac{x}{t}\big)$ is increasing (respect to $|x|$), the above estimate implies that in $[-r_0,r_0]$ the supremum of $|(te^{\frac{\cdot}{t}})'(x)|$ is less than the infimum of the absolute value of the derivative of $u^k$ in $[-2^{-k},2^{-k}]$. One can easily see that this guarantees the coincidence of functions $u^k_t$ and $te^{\frac{|x|}{t}}$ in $[-r_0,r_0]$ if we only know that 
\[te^{\frac{r_0}{t}}\leq C_k(2^{-k})^{\alpha}\,.\]  
But this follows by substituting the value of $r_0$ and using the assumption $t\leq 2^{-k+1}\,$. 

Observe then that if functions $u^k$ and correspondingly $u^k_t$ above were defined in $\rn$ (instead of $\re\,$) in the same way as in (\ref{definit}) and (\ref{definit2}), we would end up exactly with the same conclusion, thus the coincidence of $u^k_t$ and $te^{\frac{x}{t}}$ in $B(0,r_0)\,$. Let us then estimate the $p$-norm of the derivative of $u_t$:
\begin{align*}
\int_{B(0,r_0)}|Du^k_t(x)|^p\,dx\,&=\int_{B(0,r_0)}\big(e^{\frac{|x|}{t}}\big)^p\,dx\,\gtrsim 
\int_0^{r_0}z^{n-1}e^{\frac{zp}{t}}\,dz\,\\
&\gtrsim\,r_0^{n-1}\int_0^{r_0}e^{\frac{zp}{t}}\,dz\,\gtrsim r_0^{n-1}\frac{t}{p}e^{\frac{r_0p}{t}}
\\&=t^{n-1}\bigg(\log\bigg(\frac{t^{\alpha-1}C_k}{2}\bigg)\bigg)^{n-1}\frac{t}{p}\big(\frac{t^{\alpha-1}C_k}{2}\big)^p\\
&\approx t^{n+(\alpha-1)p}C_k^p\,\bigg(\log\bigg(\frac{t^{\alpha-1}C_k}{2}\bigg)\bigg)^{n-1}\,.
\end{align*}
Here (as well as below) notation $"\gtrsim"$ and $"\approx"$ mean that corresponding equality/inequality is valid up to a constant which does not depend on $t$ or $k$. 
Furthermore, 
\begin{align*}
\int_{\rn}|Du^k(x)|^p\,dx\,&\approx\int_{B(0,2^{-k})}C_k^p|x|^{(\alpha-1)p}
\approx C_k^p\int_0^{t}z^{n-1}z^{(\alpha-1)p}\,dz\,\\
&\approx t^{n+(\alpha-1)p}C_k^p\,\approx\,\frac{1}{k^{\frac{3}{2}}}. 
\end{align*}

Finally, compute that 
\begin{align*}
\bigg(\log\bigg(\frac{t^{\alpha-1}C_k}{2}\bigg)\bigg)^{n-1}&\geq
\bigg(\log\bigg(\frac{2^{(-k+1)(\alpha-1)}C_k}{2}\bigg)\bigg)^{n-1}\\
&=\bigg(\log\bigg(\frac{2^{(k-1)(1-\alpha)}2^{k(\alpha-1+\frac{n}{p})}\,}{2k^{\frac{3}{2p}}}\bigg)\bigg)^{n-1}\\
&\gtrsim\bigg(\log\bigg(\frac{2^{k\frac{n}{p}}}{k^{\frac{3}{2p}}}\bigg)\bigg)^{n-1}
\,\gtrsim k^{n-1}\,. 
\end{align*}

Combining the above estimates and substituting the value of $C_k$ we get that 
\begin{align*}
\int_{B(0,r_0)}|Du^k_t(x)|^p\,dx\,
%\gtrsim k^{n-1}\int_{\rn}|Du^k(x)|^p\,dx\,\\
\gtrsim k^{n-1}C_k^pt^{n+(\alpha-1)p}\approx k^{n-1}\frac{1}{k^{\frac{3}{2}}}\to\infty \text{ when } k\to\infty\,.
\end{align*}
Above in the last convergence we need the assumption $n\geq 3$ (this assumption is commented below). 

Summing up, we constructed a sequence of radial Sobolev-functions $u^k$ so that $(u^k)\to 0$ in $W^{1,p}(\rn)$, functions $u^k$ are constant outside the ball $B(0,2^{-k})$ and so that \text{if} $t_k$ is in range $[2^{-k},2^{-k+1}]$, then the $p$-norm of $Du^k_{t_k}$ tends arbitrarily big as $k$ tends to infinity. To verify our claim, thus the existence of single compactly supported Sobolev-function $u$ for which $\normi{Du_t}_p\to\infty$ as $t\to\infty$, we will use functions $u^k$ in the following way: Define a collection of disjoint balls $B_k=B(x_k,r_k)\subset\rn$ so that $r_k=k2^{-k}\,$ and all these balls are contained in $B(0,R)$, $R>0$. Then define function $u$ so that outside the union of these balls $u\equiv 0$ and in every $B(x_k,r_k)$ $u$ coincides with $\tilde{u}^k$, where
\[\tilde{u}^k(x)=u^k(x-x_k)-C_k2^{(-k)\alpha}\,.\]
That is, $\tilde{u}^k$ is given by translation of $u^k$ so that $\tilde{u}^k=0$ in $B(x_k,r_k)\setminus B(x_k,2^{-k})\,$.
Then it is easy to check that if $t\in [2^{-k},2^{-k+1}]$ then $u_t$ coincides in $B(x_k,r_k)$ with $u^k_t$ up to the appropriate translation. To guarantee this fact we defined $r_k=k2^{-k}$ to get enough space around every translated $u^k$. 

 Then what we computed above implies that $\normi{Du_t}_p\to\infty$ as $t\to 0\,$. 
The only thing we still have to check is that $\normi{Du}_p$ is finite. But this follows since 
\begin{equation}\label{last}
\normi{Du}_p\leq \sum_{k=1}^{\infty}\normi{Du^k}_p\leq \sum_{k=1}^{\infty}\frac{1}{k^{\frac{3}{2}}}\,<\infty\,. 
\end{equation}
\hfill $\Box$

The above theorem was proved only in the case $n\geq 3$. How about the cases $n=2$ or $n=1$? The one-dimensional case appears as a special case in many ways. Indeed, we state the following conjecture:
\begin{conjecture}
If $n=1$ and $L$ satisfies the assumptions given in the introduction, then $\normi{u'_t}_p\leq C\normi{u'}_p\,$ for any $1\leq p\leq\infty$.
\end{conjecture}
As we proved, in multidimensional case the above theorem may fail if $L$ behaves in infinity like exponential function. In the core of the argument was the fact that exponential function $L$ may 'spread' the singularity of the derivative in such a way that the measure of the set for which the bad behaviour of the gradient is copied becomes enough big to destroy the uniform $p$-integrability.
This phenomenon can not happen in one-dimensional case. However, we predict that proving the above conjecture requires a subtle argument since it seems that simple maximal-function type pointwise inequalities do not exist.

In case $n=2$ it is easy to check that above proof implies that there exists compactly supported $u\in W^{1,p}(\re^2)$ so that 
\[\limsup_{t\to 0}\normi{Du_t}_{p}\to\infty\,.\] 
However, the stronger divergence in this case seems to require 'less rough' constructions and estimates than it was used in the proof of the previous theorem. The difficulty lies in the fact that to make $\normi{Du_t}$ big it seems technically challenging to avoid of using a sequence of 'singularity functions' $u^k$ with pairwise disjoint supports as in above construction ('every scale needs its own function'). Eventually, this kind of construction requires that $\sum \normi{Du^k}_p$ is finite, say typically $\normi{Du^k}_p\lesssim\frac{1}{k}$. 
But this leads to the situation where the increment in the Sobolev-norm of the gradient, caused basically by the 'measure factor' is comparable to $k^{n-1}$ which in this case equals to $k$. In this case factors $k$ and $\frac{1}{k}$ cancel and one does not get the desired divergence. Summing up, we state the following problem:

\begin{question}
Let $1<p<\infty$ and define functions $u_t$ using proper function $L$ such that $L(x)=e^{|x|}$ when $r\geq r_0$. Does there exist compactly supported function $u\in W^{1,p}(\re^2)$ such that $\normi{Du_t}_p\to\infty$ as $t\to\infty\,$? Furthermore, one can ask if this $u$ can be even radial (this question can be posed in the case $n\geq 3$, as well).
\end{question}

\subsection*{Quasi-radiality}
We end up with a simple example relating to the necessity of the assumption (\ref{paha3a}) for our results. Suppose that $s>s'>1$, and define strictly convex function $L$ in the plane by $L(x,y)=|x|^s+|y|^{s'}\,$. Elementary calculations show that now $L$ is Legendre transform of Hamiltonian $H$, where
\begin{align*}
H(x,y)&=\bigg(\frac{s-1}{s}\bigg)\bigg(\frac{1}{s}\bigg)^{\frac{1}{s-1}}|x|^{\frac{s}{s-1}}\,+\,\bigg(\frac{s'-1}{s'}\bigg)\bigg(\frac{1}{s'}\bigg)^{\frac{1}{s'-1}}|y|^{\frac{s'}{s'-1}}\,\\
&=:C(s)|x|^{\frac{s}{s-1}}\,+\,C(s')|y|^{\frac{s'}{s'-1}}\,.
\end{align*}  
We consider this as the most natural example of the failure of (\ref{paha3a}) and in what follows we will show that if $L$ is as above then for any $2<p<\infty$ there exists compactly supported and radial function $u\in W^{1,p}(\re^2)$ such that $\normi{Du_t}_p\to \infty$ as $t\to 0$. Unless our example is in $\re^2$ it is obvious that the same phenomenon can take place also in higher dimensions.
 
Let $0<\alpha<1$ and suppose that function $u$ is compactly supported such that 
\begin{equation}
u(x,y)=|(x,y)|^{\alpha}\, \text{ if }x\in B(0,1)
\end{equation}
and $u$ is smooth outside the origin. Then, if $p>2$, $u\in W^{1,p}(\rn)$ exactly when $\alpha>\frac{p-2}{p}$. We claim that now
\[\normi{u_t}_p\to\infty\text{ as }t\to 0\,,\] 
if $\alpha$ is chosen such that $\alpha-\frac{p-2}{p}$ is positive but small enough. It turns out that it suffices to consider the set of those $(x,y)$ for which  
%The crucial point is that there exist $c_1>0$ and $c_2>0$, which especially do not depend on $t$, such that 
\begin{equation}\label{coincide}
u_t(x,y)=\frac{|x|^s}{t^{s-1}}+\frac{|y|^{s'}}{t^{s'-1}}\,. 
\end{equation}
In this set the infimum in (\ref{hopf}) is achieved at $y=0$ (or alternatively, $\frac{-x}{t}\in\mathcal{R}f(x)\,$).
We claim that there exist $c_1>0$ and $c_2>0$, which especially do not depend on $t$, such that (\ref{coincide}) holds if $(x,y)\in Q$, where  
\begin{equation}\label{crucial}
Q=\big{\{}(x,y)\in\re^2\,:\,|x|\leq c_1t^{\frac{s-1}{s-\alpha}}\,\text{ and }\,|y|\leq c_2t^{\frac{s'-1}{s'-\alpha}}\,\big{\}}. 
\end{equation}
To show this, observe that if (\ref{coincide}) does not hold, then there exists $a\in\re^2$ such that
\begin{equation}\label{alku}
|a|^{\alpha}+\frac{|a_1-x|^s}{t^{s-1}}+\frac{|a_2-y|^{s'}}{t^{s'-1}}<\frac{|x|^s}{t^{s-1}}+\frac{|y|^{s'}}{t^{s'-1}}\,.
\end{equation}
It is easy to check that we may assume that above $|a_1|\leq |x|$ and $|a_2|\leq |y|$. Then it follows from (\ref{alku}) that
\begin{align*}\label{alku2}
|a_1|^{\alpha}+|a_2|^{\alpha}\lesssim |a|^{\alpha}&<\frac{|x|^s-|x-a_1|^s}{t^{s-1}}+\frac{|y|^{s'}-|y-a_2|^{s'}}{t^{s'-1}}\\
&\lesssim \frac{|a_1||x|^{s-1}}{t^{s-1}}+\frac{|a_2||y|^{s'-1}}{t^{s'-1}}\,,
\end{align*}
implying that either
\begin{align*}
&|a_1|^{1-\alpha}\gtrsim \frac{t^{s-1}}{|x|^{s-1}}\, \text{ or }\,|a_2|^{1-\alpha}\gtrsim \frac{t^{s'-1}}{|y|^{s'-1}}\,.\\
\end{align*}
Finally, using again facts $|a_1|\leq |x|$ and $|a_2|\leq |y|$ we get that either
\begin{align*}
&|x|\gtrsim t^{\frac{s-1}{s-\alpha}}\, \text{ or }\,|y|\gtrsim t^{\frac{s-1}{s-\alpha}}\,.  
\end{align*}
This guarantees that (\ref{coincide}) holds in $Q$ for constants $c_1,c_2>0$, independent on $t$. Then, one can estimate $\normi{Du_t}_{p}$ from below simply by 
\begin{align*}
&\int_{\re^2}|Du_t(x,y)|^p\,d(x,y)\,\geq\int_{Q}|D_xu_t(x,y)|^p\,d(x,y)\,\\
\geq & \int_{0}^{c_1t^{\frac{s-1}{s-\alpha}}}\int_{0}^{c_2t^{\frac{s'-1}{s'-\alpha}}}s^p\frac{|x|^{(s-1)p}}{t^{(s-1)p}}\,dy\,dx
\gtrsim \int_{0}^{c_1t^{\frac{s-1}{s-\alpha}}}t^{\frac{s'-1}{s'-\alpha}}\frac{|x|^{(s-1)p}}{t^{(s-1)p}}\,dx\,\\
\gtrsim & \,t^{\frac{s-1}{s-\alpha}}t^{\frac{s'-1}{s'-\alpha}}\frac{t^{\frac{s-1}{s-\alpha}(s-1)p}}{t^{(s-1)p}}
=t^{\frac{s-1}{s-\alpha}+\frac{s'-1}{s'-\alpha}+\frac{p(s-1)(\alpha-1)}{s-\alpha}}=:t^{\beta}\,.
\end{align*}
By elementary calculations, $\beta<0$ is equivalent with  
\begin{equation}\label{joo}
(s'-1)(s-\alpha)+(s-1)(s'-\alpha)<p(s-1)(1-\alpha)(s'-\alpha)\,.
\end{equation}
We are ready if we only can show that $\beta<0$ for suitable $1>\alpha>\frac{p-2}{p}$. This in turn follows if 
(\ref{joo}) holds with substitution $\alpha=\frac{p-2}{p}$. Then $1-\alpha=\frac{2}{p}$ and (\ref{joo}) is equivalent with
\begin{align*}\label{joo2}
&(s'-1)(s-\alpha)+(s-1)(s'-\alpha)<2(s-1)(s'-\alpha)\,\\
\Longleftrightarrow\,\,&(s'-1)(s-\alpha)<(s-1)(s'-\alpha)\,\\
\Longleftrightarrow\,\,&-\alpha s'-s<-\alpha s-s'\,\\
\Longleftrightarrow\,\,&\alpha =\frac{p-2}{p}<1\,.
\end{align*}
This verifies our claim.

\end{document}